\documentclass[12pt]{article}
\usepackage[english]{babel}
\usepackage{amssymb}
\usepackage{inputenc}
\begin{document}
\begin{center}
\textsc{On \v{C}ech-completeness of the space of order-preserving
functionals}

Sh. A. Ayupov\footnote{Institute of Mathematics, National
university of Uzbekistan, Tashkent; and the Abdus Salam Centre for
Theoretical Physics (ICTP), Trieste, Italy}, A. A.
Zaitov\footnote{Institute of Mathematics, Uzbekistan, Tashkent;
and THGACSS}

\end{center}

\begin{abstract}

In this paper we establish that if a  Tychonoff space $X$ is
\v{C}ech-complete then the space $O_\tau(X)$ of all $\tau$-smooth
order-preserving, weakly additive and normed functionals is also
\v{C}ech-complete .

\end{abstract}



\section{Introdiction}

Nowadays the study of nonlinear functionals becomes one of
attractive fields in topology and functional analysis. For
example, in the papers [1] -- [3] so called \textit{max-plus} or
idempotent linear functionals have been considered. Such
functionals have  various applications in mathematical physics,
economics, dynamical optimization and other fields. Another
example of nonlinear functionals is given  by the notion of
capacity on compacta introduced by Choquet [4]. In [5] various
classes of nonlinear functionals have been considered and
investigated . Order-preserving, weakly additive functionals on
the algebra of continuous functions on a given compact Hausdorff
space were considered in [6]. Note that nonnegative linear
functionals, \textit{max-plus} functionals and capacities are
order-preserving and weakly additive.

In this paper we consider the functor $O_\tau$ of order-preserving
$\tau$-smooth functionals in the category $Tych$ of Tychonoff
spaces and their continuous mappings, and for a given Tychonoff
space $X$ we establish that if $X$ is \v{C}ech-complete, then the
space $O_\tau(X)$ of $\tau$-smooth order-preserving functionals is
also \v{C}ech-complete .

Earlier proven the Baire category theorem  for the class of all
completely metrizable spaces has been generalized for the case of
\v{C}ech-complete spaces. Therefore one of the important questions
in the theory of covariant functors is whether a functor preserves
\v{C}ech-completeness of spaces.

We denote by $\mathbb{R}$ the set of all real numbers, by
$\mathbb{N}$ the set of all nonnegative integers,  spaces mean
topological spaces, $\overline{A}$ is the closure of $A$ in the
considering space, compact is compact Hausdorff topological space,
maps (or mappings) are continuous ones.

\section{On properties of order-preserving functionals }

Let $X$ be a compact and let $C(X)$ be the Banach algebra of all
continuous real-valued functions with the usual algebraic
operations and with the $\sup$-norm. For functions $\varphi,\
\psi\in C(X)$ we shall write $\varphi\leq\psi$ if
$\varphi(x)\leq\psi(x)$ for all $x\in X.$ If $c\in \mathbb{R}$
then by $c_X$ we denote the constant function which is identically
equal to $c.$

We recall that a functional $\mu : C(X)\rightarrow \mathbb{R}$ is
said [6] to be:

1) \textit{order-preserving} if for any pair $\varphi, \psi\in
C(X)$ of functions the inequality $\varphi\leq\psi$ implies
$\mu(\varphi)\leq\mu(\psi)$;

2) \textit{weakly additive} if
$\mu(\varphi+c_X)=\mu(\varphi)+c\mu(1_X)$ for all $\varphi\in
C(X)$ and $c\in R$;

3) \textit{normed} if $\mu(1_X)=1$.

For a compact $X$ we denote by $W(X)$ the set of all
order-preserving weakly additive functionals $\mu :
C(X)\rightarrow \mathbb{R}$, and by $O(X)$ we denote the set of
all normed functionals from $W(X)$. It is clearly that $W(X)$ is a
subset in the product $\mathbb{R}^{C(X)}$. Equip $W(X)$ with the
induced topology, which coincides with the pointwise convergence
topology on $W(X)$. Note that in this topology the space $O(X)$ is
compact, for any given compact $X$.

For a given continuous map $f:X\rightarrow Y$ between compacta by
$$
O(f)(\mu)(\varphi)=\mu(\varphi\circ f),\ \qquad \varphi\in C(Y),
$$
we define the map $O(f):O(X)\rightarrow O(Y)$. Thus we get the
functor $O$ acting in the category $Comp$.

Order-preserving, weakly additive and normed functionals for
brevity we will call as order-preserving functionals.

Note that according to proposition 1 in [7] each order-preserving
functional is continuous.

Let $A$ be a closed subset of a given compact $X$. Since the
functor $O:Comp\rightarrow Comp$ is monomorphic [6], i. e.
preserves embeddings of closed subsets of compacta, we may assume
that $O(A)$ is a subset of the compact $O(X).$ Moreover the
embedding of $O(A)$ into $O(X)$ is given by means of the mapping
$O(i_A): O(A)\rightarrow O(X),$ where $i_A: A\rightarrow X$ is the
identical embedding. Therefore the following notion is well
defined.

Given a closed subset $A$ of the compact $X$, an order-preserving
functional $\mu\in O(X)$ is said to be \textit{supported} on $A$
if $\mu\in O(A)$. The set
$$
\mbox {supp}\mu=\cap\{A:\ \mu\in O(A) \mbox{ and } A \mbox{ is
closed in } X\}
$$
is called the \textit{support} of the order-preserving functional
$\mu.$

Let $X$ be a compact. For a closed $F\subset X$, an open $U\subset
X$, real numbers $\alpha$ and $c$ put
$$\begin{array}{ll}
\mbox{if } \alpha\geq 0, \mbox{ then }
\mu(\alpha\chi_F+c_X)=\inf\{\mu(\varphi): \varphi\in C(X),\
\varphi\geq\alpha\chi_F\}+c\mu(1_X), &(1^+)\\
\mbox{if }  \alpha< 0, \mbox{ then }
\mu(\alpha\chi_F+c_X)=\sup\{\mu(\varphi): \varphi\in C(X),\
\varphi\leq\alpha\chi_F\}+c\mu(1_X), &(1^-)\\
\mbox{if } \alpha\geq 0, \mbox{ then }
\mu(\alpha\chi_U+c_X)=\sup\{\mu(\varphi): \varphi\in C(X),\
\varphi\leq\alpha\chi_U\}+c\mu(1_X), &(2^+)\\
\mbox{if } \alpha< 0, \mbox{ then }
\mu(\alpha\chi_U+c_X)=\inf\{\mu(\varphi): \varphi\in C(X),\
\varphi\geq\alpha\chi_U\}+c\mu(1_X). &(2^-)
\end{array}$$
In this manner  the functional $\mu$  given on $C(X)$ can be
extended onto the following set
$$
\mathcal{X}(X)=C(X)\cup\{\alpha\chi_A+c_X: \alpha, c\in
\mathbb{R}, A \mbox{ is open or closed in } X\}.
$$
Obviously the set $\mathcal{X}(X)$ is an $A$-\textit{subspace} of
the space $B(X)$ of all bounded functions on $X$, i. e. $0_X\in
\mathcal{X}(X)$ and for any $\zeta\in \mathcal{X}(X)$ and each
$c\in \mathbb{R}$ we have $\zeta+c_X\in \mathcal{X}(X)$. Hence by
 the version of the Hahn-Banach theorem [9] it follows that $\mu$
 has an extension $\widetilde{\mu}$ on $B(X)$ such that
$\widetilde{\mu}|\mathcal{X}(X)=\mu$.

The following result has been proved in [8, Th. 1.3].

\textbf{Theorem 1.} \textsl{Let $A$ be an arbitrary subspace of
the given compact $X$, $\mu\in O(X)$. Then for each extension
$\widetilde{\mu}$ satisfying $(1^{\pm})-(2^{\pm})$, the following
conditions are equivalent:}

$(i)$ $\widetilde{\mu}(n\chi_A)=0$ \textsl{for all integers} $n$.

$(ii)$ $\widetilde{\mu}(\alpha\chi_A)=0$ \textsl{for all}
$\alpha\in \mathbb{R}$.

$(iii)$ $\widetilde{\mu}(\alpha\chi_{X\setminus A})=\alpha$
\textsl{for all} $\alpha\in \mathbb{R}$.

$(iv)$ $\widetilde{\mu}(\varphi\chi_A)=0$ \textsl{for all}
$\varphi\in C(X)$.

$(v)$ $\widetilde{\mu}(\varphi\chi_{X\setminus A})=\mu(\varphi)$
\textsl{for all} $\varphi\in C(X)$.

$(vi)$ $\mbox{supp}\mu\subset \overline{X\setminus A}$.

We need the following statement.

\textbf{Proposition 1.}  \textsl{For each open subset $U$ of a
given compact $X$ and for an arbitrary $\varepsilon>0$ the sets}
$$
\langle U; \varepsilon\rangle=\{\mu\in O(X):\
\mu(\chi_U)>\varepsilon\}\  \mbox{ and }  \ \langle U;
-\varepsilon\rangle=\{\mu\in O(X):\ \mu(-\chi_U)<-\varepsilon\}
$$
\textsl{are open in $O(X)$ with respect to pointwise convergence
topology.}

\textsc{Proof.}  For each $\varphi\in C(X)$ and $\varepsilon>0$
put
$$
\langle\varphi;\ \varepsilon\rangle=\{\mu\in O(X):\
\mu(\varphi)>\varepsilon\}.
$$
From the definition of pointwise convergence topology follows that
the set $\langle\varphi;\ \varepsilon\rangle$ is open in $O(X)$.

Consider the following set
$$
\Phi=\{\varphi\in C(X):\ 0_X \leq \varphi\leq\chi_U,\ \mbox{ there
exists } \mu\in O(X) \mbox{ such that }
\mu(\varphi)>\varepsilon\}.
$$
It is clear that $\Phi\neq\varnothing$.

We have
$$
\langle U; \varepsilon\rangle=\bigcup_{\varphi\in \Phi
}\langle\varphi;\ \varepsilon\rangle,
$$
i. e. the set $\langle U; \varepsilon\rangle$ is open in $O(X)$.
Analogously, the set $\langle U; -\varepsilon\rangle$ is also open
in $O(X)$. Proposition 1 is proved.

We note that the sets
$$ \langle U;
\varepsilon\rangle_{\lambda}=\{\mu\in O(X):\
\mu(\lambda\chi_U)>\varepsilon\},\  \qquad \lambda>0,
$$
and
$$
\langle U; -\varepsilon\rangle_{\lambda}=\{\mu\in O(X):\
\mu(\lambda\chi_U)<-\varepsilon\}, \qquad \lambda<0,
$$
are also open in $O(X)$.

\section{Order-preserving $\tau$-smooth functionals}

Let $X$ be a Tychonoff space and let $C_b(X)$ be the algebra of
all bounded continuous real-valued functions with the pointwise
algebraic operations. For a function $\varphi\in C_b(X)$ put
$\|\varphi\|=\sup\{|\varphi(x)|:\ x\in X\}.$ $C_b(X)$ with this
norm is a Banach algebra. For a net $\{\varphi_\alpha\}\subset
C_b(X)$ the notation $\varphi_\alpha \downarrow 0_X$ means that
for every point $x\in X$ at $\beta \succ \alpha$  implies the
inequality $\varphi_\alpha(x) \geq \varphi_\beta(x)$, and
$\lim\limits_\alpha \varphi_\alpha(x)=0_X.$ In this case we say
that $\{\varphi_\alpha\}$ is a \textit{monotone decreasing net
pointwise convergent to zero}.

For a Tychonoff space $X$ we denote by $\beta X$  its
Stone-\v{C}ech compact extension. Given any function $\varphi\in
C_b(X)$ consider its continuous extension $\widetilde{\varphi}\in
C(\beta X).$ This gives an isomorphism between the spaces $C_b(X)$
and $C(\beta X)$. Moreover, $\|\widetilde{\varphi}\|=\|\varphi\|,$
i. e. this isomorphism is an isometry, and the topological
properties of the above spaces coincide. Therefore one may
consider any function from $C_b(X)$ as an element of $C(\beta X).$

\textbf{Definition 1.} An order-preserving functional $\mu\in
W(\beta X)$ is said to be \textit{strongly $\tau$-smooth} if
$\mu(\pm\varphi _\alpha)$ $\rightarrow 0$ for each monotone  net
$\{\varphi_\alpha\} \subset C(\beta X)$ decreasing to zero on $X.$

For brevity strongly $\tau$-smooth functionals we will call
$\tau$-smooth functionals.

For a Tychonoff space $X$ we denote by $W_\tau(X)$ the set of all
$\tau$-smooth order-preserving functionals from $W(\beta X).$  The
set $W_\tau(X)$ is equipped with the pointwise convergence
topology. The base of neighborhoods of a functional $\mu \in
W_\tau (X)$ in the pointwise convergence topology consists of the
sets
$$
\langle \mu; \varphi_1,...,\varphi_k; \varepsilon \rangle = \{\nu
\in W(\beta X): |\nu(\varphi_i) - \mu(\varphi_i)| < \varepsilon \}
\cap W_\tau (X)
$$
where $\varphi_i \in C(\beta X),$ $i=1,...,k$ and $\varepsilon >
0$.

Put
$$
O_\tau(X)=\{\mu\in W_\tau(X): \mu(1_X)=1\},
$$
It is easy to see that equipped with pointwise convergence
topology $O_\tau(X)$ is Tychonoff space.

Let $X$ and $Y$ be Tychonoff spaces and let $f:X\rightarrow Y$ be
a continuous map. Put
$$
O_\tau(f)=O(\beta f)|O_\tau(X),
$$
where $\beta f: \beta X\rightarrow\beta Y$ is the Stone-\v{C}ech
extension of $f.$

Note that the map $O_\tau(f):O_\tau(X)\rightarrow O_\tau(Y)$ is
defined correctly and it is also continuous. Thus we obtain the
operation $O_\tau$ as functor acting in the category $Tych$ of
Tychonoff spaces and their maps [8].

The following theorem, which gives an equivalent definition of
$\tau$-smoothness of order-preserving functionals, was proved in
[8, Th. 2.2].

\textbf{Theorem 2.} \textsl{A functional $\mu\in O(\beta X)$ is
$\tau$-smooth if and only if
$\mu(\varphi_{\alpha})\rightarrow\mu(\varphi)$ and
$\mu(-\varphi_{\alpha})\rightarrow\mu(-\varphi)$ for any monotone
increasing net of continuous functions
$\{\varphi_{\alpha}\}\subset C(\beta X)$, pointwise converging to
$\varphi\in C(\beta X)$ on $X$}.

Now we shall prove a result which can be consider as another
equivalent definition of $\tau$-smoothness of order-preserving
functionals.

\textbf{Theorem 3.} \textsl{Functional $\mu\in O(\beta X)$ is
$\tau$-smooth if and only if $\mu(\lambda\chi_K)=0$ for each
compact $K\subset \beta X\setminus X$ and for every $\lambda\in
\mathbb{R}$.}

\textsc{Proof}. Let $\mu\in O(\beta X)$ be an arbitrary
$\tau$-smooth functional and let $K\subset \beta X\setminus X$ be
a compact. Consider a monotone decreasing net
$\Phi=\{\varphi_{\alpha}\}\subset C(\beta X)$ such that
$\varphi_{\alpha}\geq \lambda\chi_K$ and
$\inf\{\mu(\varphi_{\alpha}): \varphi_{\alpha}\in
\Phi\}=\inf\{\mu(\varphi): \varphi\in C(\beta X),\ \varphi\geq
\lambda\chi_K\}\equiv\mu(\lambda\chi_K)$, where $\lambda>0$. Then
$\mu(\varphi_{\alpha})\rightarrow 0$.  Hence,
$\mu(\lambda\chi_K)=0$. From the construction immediately follows
that $\{\mu(-\varphi_{\alpha})\}$ is an increasing net and
$\mu(-\lambda\chi_K)=\sup\{\mu(-\varphi_{\alpha}):
\varphi_{\alpha}\in \Phi\}$. By the definition
$\mu(-\varphi_{\alpha})\rightarrow 0$. But then it should be
$\mu(-\lambda\chi_K)=0$.

Thus for any compact $K\subset \beta X\setminus X$ and for all
$\lambda\in \mathbb{R}$ the equality $\mu(\lambda\chi_K)=0$ takes
place.

Conversely, let $\mu\in O(\beta X)$ be  a functional such that
$\mu(\lambda\chi_K)=0$ for any compact $K\subset \beta X\setminus
X$ and for an arbitrary $\lambda\in \mathbb{R}$. Consider an
arbitrary net $\Phi=\{\varphi_{\alpha}\}\subset C(\beta X)$,
$\varphi_{\alpha}\downarrow 0_X$. It is sufficient to check the
value $\lambda=1$. Without loss of generality we may assume that
$\varphi_{\alpha}\leq 1_{\beta X}$ for all $\alpha$.

For an arbitrary real number $\varepsilon>0$ consider the set
$$Z_{\alpha}=\{x\in \beta X: \varphi_{\alpha}(x)\geq \varepsilon\}.$$
It is clear that
$K_{\varepsilon}:=\bigcap\limits_{\alpha}Z_{\alpha}\subset \beta
X\setminus X$. Moreover it is easy to see that for $\alpha\prec
\beta$ the inclusion $Z_{\alpha}\supset Z_{\beta}$ is valid.
Therefore $\chi_{Z_{\alpha}}\geq\chi_{Z_{\beta}}$ for $\alpha\prec
\beta$. Hence, $\mu(\chi_{Z_{\alpha}})\geq\mu(\chi_{Z_{\beta}})$.
In other words the net $\{\mu(\chi_{Z_{\alpha}})\}$ is decreasing.
From the equality
$K_{\varepsilon}=\bigcap\limits_{\alpha}Z_{\alpha}$ it follows
that the net $\{\mu(\chi_{Z_{\alpha}})\}$ is  bounded from below
by $\mu(\chi_{K_{\varepsilon}})$. Moreover,
$\mu(\chi_{Z_{\alpha}})\rightarrow\mu(\chi_{K_{\varepsilon}})$.
Hence one has
$$\mu(\varphi_{\alpha})=\widetilde{\mu}(\varphi_{\alpha})=
\widetilde{\mu}(\varphi_{\alpha}\chi_{Z_{\alpha}}+
\varphi_{\alpha}\chi_{\beta X\setminus
Z_{\alpha}})\leq(\mbox{since }  \widetilde{\mu} \mbox{ is
order-preserving})\leq$$
$$\leq\widetilde{\mu}(\chi_{Z_{\alpha}}+\varepsilon_{\beta X})=
(\mbox{by the weak additivity of }
\widetilde{\mu})=\widetilde{\mu}(\chi_{Z_{\alpha}})+\varepsilon\mu(1_{\beta
X})=$$
$$=\mbox{(by definition of } \widetilde{\mu}) =\mu(\chi_{Z_{\alpha}})+\varepsilon\mu(1_{\beta X})\rightarrow
\mu(\chi_{K_{\varepsilon}})+ \varepsilon\mu(1_{\beta
X})=\varepsilon\mu(1_{\beta X}).$$

Now taking $\varepsilon$ convergent to zero we have
$\mu(\varphi_{\alpha})\rightarrow 0$.

One can similarly  establish the convergence
$\mu(-\varphi_{\alpha})\rightarrow 0$. Theorem 3 is proved.

According to Theorems 1 and 3 we can write
$$
\begin{array}{l}
O_\tau(X)=\{\mu\in O(\beta X):\ \mu(n\chi_K)=0  \mbox{ for each
compact }\ K\subset \beta X \mbox{ such }\\
\hspace{4cm}\mbox{ that }\ K\cap X=\varnothing, \ \mbox{ and for
all integers } \ n\}. \ \qquad (3)
\end{array}
$$

\section{The main result}

Recall that a topological space $X$ is called to be
\textit{\v{C}ech-complete} if $X$ is a Tychonoff space and for
some (or, equivalently, for every) compactification $bX$ of the
space $X$ the remainder $bX\setminus X$ is an $F_\sigma$-set in
$bX$ [10, P. 196].

The following statement is the main result of the paper.

\textbf{Theorem 4.} \textsl{If a Tychonoff space $X$ is
\v{C}ech-complete then the space $O_\tau(X)$ is also
\v{C}ech-complete.}

\textsc{Proof.} Let a  Tychonoff space $X$ be \v{C}ech-complete.
Then the remainder $\beta X\setminus X$ is an $F_\sigma$-set in
$\beta X$, i. e. in $\beta X$ there exist closed subsets $K_{n}$,
$n\in \mathbb{N}$ such that $\beta X \setminus
X=\cup_{n=1}^{\infty} K_{n}$. Hence  $X=\cap_{n=1}^{\infty}U_{n}$
is a $G_{\delta}$-set in $\beta X$, where $U_{n}=\beta X\setminus
K_n$ are open in $\beta X$ sets.

Put $F_n=\bigcup\limits_{i=1}^{n}K_n$, $V_n=\beta X\setminus F_n$.
Then $F_n\subset F_{n+1}$ for all $n\in \mathbb{N}$ and $\beta
X\setminus X=\cup_{n=1}^{\infty}F_n$.

According to (3) for each $\mu\in O_\tau(X)$ we have
$\mu(\lambda\chi_{F_n})=0$ for all $n=1,2,...$ and $\lambda\in
\mathbb{N}$. From this  by Theorem 1 it follows that
$\mu(\lambda\chi_{V_n})=\lambda$. Consider decreasing sequence
$\{r_{n}\}$ of positive numbers such that
$\lim\limits_{n\rightarrow\infty}r_{n}=0$. Then for any
$\lambda\in \mathbb{N}$ we have $\mu\in \langle V_n;\
(\lambda-r_{n})\rangle_{\lambda}\cap\langle V_n;\
-(\lambda-r_{n})\rangle_{\lambda}$ (here without loss of
generality we assumed that $r_{n}<\lambda$ for all
$n\in\mathbb{N}$). From this we have the following inclusion
$$
O_\tau(X)\subset\bigcap\limits_{\lambda=1}^{\infty}\bigcap\limits_{n=1}^{\infty}(\langle
V_n;\ (\lambda-r_{n})\rangle_{\lambda}\cap\langle V_n;\
-(\lambda-r_{n})\rangle_{\lambda}).\eqno(4)
$$

Let now $\mu\in
\bigcap\limits_{\lambda=1}^{\infty}\bigcap\limits_{n=1}^{\infty}(\langle
V_n;\ (\lambda-r_{n})\rangle_{\lambda}\cap\langle V_n;\
-(\lambda-r_{n})\rangle_{\lambda})$. Then for all positive
integers $\lambda$ and $n$ we will have
$\mu(\lambda\chi_{V_n})>\lambda- r_{n}$ and
$\mu(-\lambda\chi_{V_n})<-\lambda+ r_{n}$. Hence
$\mu(\lambda\chi_{\bigcap\limits_{n=1}^{\infty}V_n})\geq\lambda$
and
$\mu(-\lambda\chi_{\bigcap\limits_{n=1}^{\infty}V_n})\leq-\lambda$,
because $\{V_n\}$ is a decreasing (with respect to inclusion)
sequence. But then
$\mu(\lambda\chi_{\bigcap\limits_{n=1}^{\infty}V_n})=\lambda$ for
all integers $\lambda$, since $\mu$ is order-preserving and
normed. Hence
$\mu(\lambda\chi_{\bigcup\limits_{n=1}^{\infty}F_n})=0$.

If $K\subset \beta X\setminus X$ is an arbitrary compact then
$\mu(\lambda\chi_{K})=0$ for all integers $\lambda$, since
$K\subset\bigcup\limits_{n=1}^{\infty}F_n$. Therefore (3) implies
that $\mu\in O_{\tau}(X)$. This means that
$$
\bigcap\limits_{\lambda=1}^{\infty}\bigcap\limits_{n=1}^{\infty}(\langle
V_n;\ (\lambda-r_{n})\rangle_{\lambda}\cap\langle V_n;\
-(\lambda-r_{n})\rangle_{\lambda})\subset O_\tau(X).\eqno(5)
$$

The inclusions (4) and (5) imply the equality
$$
O_\tau(X)=\bigcap\limits_{\lambda=1}^{\infty}\bigcap\limits_{n=1}^{\infty}(\langle
V_n;\ (\lambda-r_{n})\rangle_{\lambda}\cap\langle V_n;\
-(\lambda-r_{n})\rangle_{\lambda}).\eqno(6)
$$

The equality (6) means that $O_\tau(X)$ is a $G_\delta$-set in
$O(\beta X)$. Therefore the remainder $O(\beta X)\setminus
O_\tau(X)$ is an $F_\sigma$-set in the compactification
$bO_\tau(X)$ defined as $bO_\tau(X)=O(\beta X)$. Theorem 4 is
proved.

\end{document}